\documentstyle[12pt]{article}

\font\mathbb=msbm10 scaled \magstephalf

\def\Z{\hbox{\mathbb Z}}

\def\R{\hbox{\mathbb R}}
\def\CA{{\cal A}}
\def\CB{{\cal B}}

\def\qfact#1{\hbox{$#1 !_q$}}
\def\qbino#1#2{\left[ \begin{array}{c} #1 \\ #2 \end{array} \right]_q}
\def\beq{\begin{equation}}
\def\eeq{\end{equation}}
\newcommand{\bea}{\begin{eqnarray}}
\newcommand{\eea}{\end{eqnarray}}

\newtheorem{obs}{Observation}
\newtheorem{lem}{Lemma}

\def\zt{\Z_N}

\title{On the tensor product construction for $q$-differential algebras}
\author{
Andrzej Sitarz\thanks{e-mail:sitarz@higgs.physik.uni-mainz.de} \\ \ \\
{\em Institut f\"ur Physik, Johannes-Gutenberg Universit\"at} \\
{\em  55099 Mainz, Germany }
}

\begin{document}
\parindent 0pt

\maketitle

{\bf Abstract:}
We show that for $q\not=-1$ the $q$-graded tensor product fails to
preserve the $q$-differential structure of the product algebra and
therefore there is no {\em natural} tensor product construction for 
$q$-differential algebras.


\section{Introduction}

Recently there has been much interest in noncommutative
generalizations of the {\em standard} differential geometry.
An interesting point is the attempt to study the $\zt$-graded
differential structure and the $N$-nilpotent ($d^N \equiv 0$)
linear maps satisfying a deformed Leibniz rule 
\cite{MDV,MK2}. This could be extended to a $q$-generalization 
of complexes and cohomology \cite{MK}. 

Special examples of ${\Z}_3$-graded structures have been 
discussed earlier (see \cite{RK} and references therein)
with many examples and a view to possible physical 
realizations \cite{RK2,KRA}.

Although the generalization of differential algebras opens
new possibilities, one should not forget that it also brings  
difficulties. Not every construction that is possible in
the classical geometry can be automatically extended to a more
general setup. This Letter discusses one of the significant 
constructions, the tensor product of $q$-differential
algebras, proving that within the commonly assumed rules 
there is no room for it.

\section{$q$-differential algebras}

A $\Z$ graded $q$-differential algebra $\Omega_q$ is an
unital algebra with the following properties:
\begin{itemize}
\item $\Omega = \bigoplus_{i=0}^{\infty} \Omega^i$ and 
$\Omega^i \cdot \Omega^j \subset \Omega^{i+j}$ ($\Z$-grading),
\item there exists a linear operator $d: \Omega^i \to \Omega^{i+1}$ such
that $d^N \equiv0$,
\item $d$ obeys a $q$-graded Leibniz rule, for $q \not=1$ and $q^N=1$:
$$ d (\omega^i \omega^j) = (d \omega^i) \omega^j + q^i \omega^i 
(d \omega^j)$$
where $\omega^{i,j} \in \Omega^{i,j}$.
\end{itemize}

Note that $q$ can be {\em any primitive $N$-th root of unity}, i.e.,
$N$ is the smallest integer such that $q^N=1$. This follows from
the compatibility of the $q$-Leibniz rule with
the nilpotency of $d$. Iterating the Leibniz rule we obtain:

\begin{equation}
d^N(ab) = \sum_{n=0}^{N} \qbino{N}{n} d^n a \, d^{N-n} b, 
\label{lrq}\end{equation}
where

$$ \qbino{N}{n} = { \qfact{N} \over \qfact{n} \qfact{(N-n)} },$$

and 

$$ \qfact{n} = \frac{1-q^n}{1-q} \cdots \frac{1-q^2}{1-q}. $$

It is easy to see that for $q$ a primitive root of unity the $q$-deformed
binomial vanishes for $0<n<N$ and therefore (\ref{lrq}) agrees with
$d^N=0$.

Of course, for $N=2$ ($q=-1$ being the only possibility), 
we recover the standard differential algebras. 
However, for $N>2$ we have a whole family of differential
algebras, which differ by the value of $q$ in the Leibniz rule, as
for different values of $q$ (for the same $N$) we have different objects.

We can make this statement more precise using the notion of homomorphisms 
between $q$-differential algebra, which are algebra homomorphisms
preserving the grading and commuting with $d$:

\begin{obs}
Let $\Omega_{q}$ and $\Omega_p$ be two $q$- and $p$-differential
algebras for $q^N=p^N=1$. If there exists a differential homomorphism
$\psi: \Omega_q \to \Omega_p$ then $p=q$.
\end{obs}

The proof of this observation is simple, one should apply $\phi$ to
the Leibniz rule in $\Omega_q$ and compare the result with the
Leibniz rule in $\Omega_p$.

If the subalgebra of $0$-forms has a star operation then it extends
naturally (provided that $\Omega^n$ is generated as a bimodule by 
$\Omega^{n-1}$) on the whole $\Omega$ by the rule:

$$ (d \omega^i)^\star = q^{-i} d (\omega^\star). $$

\section{Tensor products of $q$-differential algebras}

So far, we have seen that the definition and properties of
the standard differential algebras extend without problems to
the $q$-differential case. Our task will now be to verify whether
the tensor product construction extends to the $q$-differential
case.

The construction of the differential algebra over the tensor
product of algebras (which corresponds to the construction
of differential geometry over Cartesian products of manifolds) 
uses the graded tensor product.

Let us remind that for the ordinary differential algebras $\Omega(\CA)$
and $\Omega(\CB)$ we construct $\Omega(\CA \tilde{\otimes} \CB)$ in the
following way. We take $\Omega(\CA) \otimes \Omega(\CB)$ as a vector
space, with the following product:
$$ (\omega_1 \otimes \eta_1) (\omega_2 \otimes \eta_2) =
   (-1)^{|\eta_1||\omega_2|} \left( \omega_1\omega_2 \otimes \eta_1 \eta_2
\right)
$$
for $\omega_{1,2} \in \Omega(\CA)$ and $\eta_{1,2} \in \Omega(\CB)$

Then one finds that the following linear operator $d$: 
\begin{equation}
d (\omega \otimes \eta) = (d\omega \otimes \eta) + (-1)^{|\omega|} 
(\omega \otimes d\eta), 
\label{defd}
\end{equation}
satisfies both the graded Leibniz rule and $d^2=0$, and therefore our
graded tensor product of differential algebras is a differential
algebra itself.

In this section we shall attempt to generalize the above scheme for
$q$-differential algebras.

\subsection{$q$-graded tensor product}

First, we have to introduce the notion of a braided tensor product, 
with braiding set by a root of unity $q$, $q^N=1$. 

Let us assume that ${\cal A}$ and  ${\cal B}$ are $\Z$-graded algebras.
Then we introduce the tensor product ${\cal A} \otimes {\cal B}$ 
with the following algebra structure:

$$ ( a_1 \otimes b_1)(a_2 \otimes b_2) = q^{|b_1||a_2|} (a_1 a_2 \otimes
b_1 b_2), $$

where $q$ is any $N$-th  root of unity, and $|a|$ denotes the grades
of the algebra elements. We shall denote this tensor product at $\otimes_q$.
For the proof that the construction is well-defined and the above defined
braided tensor multiplication $\otimes_q$ is associative we refer to
\cite{Majid} (In the case of the above braiding one uses often the name
of {\em anyonic vector spaces} and {\em anyonic tensor product}).  

Before we proceed with the extension of differential structures on the
tensor product, let us make a simple observation:

\begin{lem} ${\cal A} \otimes_q {\cal B}$ is isomorphic to
${\cal B} \otimes_{\bar{q}} {\cal A}$.
\end{lem}

{\bf Proof:} Let us define $\phi(a \otimes_q b) = 
q^{-|a||b|} b \otimes_{\bar{q}} a$.

Then
\begin{eqnarray*}
\phi\left( (a_0 \otimes_q b_0)(a_1 \otimes_q b_1) \right) = \\
\phi \left( q^{|b_0||a_1|} a_0 a_1 \otimes_q b_0 b_1 \right) =\\
q^{-(|a_0|+|a_1|)(|b_0|+|b_1|)} q^{|b_0||a_1|} 
\left( b_0 b_1 \otimes_{\bar{q}} a_0 a_1 \right) = \\
q^{-|a_0||b_0| - |a_1||b_1|} \left( (b_0 \otimes_{\bar{q}} a_0)
 (b_1 \otimes_{\bar{q}} a_1) \right) = \\
\phi( a_0 \otimes_q b_0) \phi(a_1 \otimes_q b_1).
\end{eqnarray*}

Therefore the anyonic tensor product (for $q \not= -1$) is
strictly noncommutative (for $q=-1$ it is $\Z_2$-graded
commutative).

Now we shall prove the main lemma.

\begin{lem}
Let $\CA$ and $\CB$ be $\Z$-graded $q$-differential algebras, with
the deformed Leibniz rule determined by $q_A$ and $q_B$ respectively,
both primitive $N$-th roots of unity. Then, no anyonic tensor
product admits a $q$-differential structure.
\end{lem}

To prove this {\em no-go} lemma we shall verify the restrictions
set upon the differential structure on $\CA \otimes_q \CB$, 
where $q$ is again a $N$-th root of unity. 

Let us rewrite the Leibniz rules for $\CA$ and $\CB$:

\begin{eqnarray}
d( a_1 a_2 ) & = & da_1 a_2 + q_A^{|a_1|} a_1 da_2, \;\;\;\;
a_1,a_2 \in \Omega({\cal A}), \\
d( b_1 b_2 ) & = & db_1 b_2 + q_B^{|b_1|} b_1 db_2, \;\;\;\;
b_1,b_2 \in \Omega({\cal B}).
\end{eqnarray}

Note, that a'priori $q_A$ might be different from $q_B$.

Now, we pose the following question: can one find a $q$, $p$, $s$
(all $N$-th roots of unity) such that the $q$-tensor product of
$\CA$ and $\CB$, admits the differential structure defined as:

$$ d (a \otimes_q b) = da \otimes_q b + s^{|a|} a \otimes_q db $$

and this $d$ obeys the $p$-deformed Leibniz rule.

We leave as much parameters free as possible, as a'priori we know
nothing about the restrictions (if any) that may appear between
them.

Of course, from the natural embedding:
\begin{eqnarray}
{\cal A} & \sim  \CA \otimes_q 1 & \hookrightarrow  {\cal A} \otimes_q {\cal B} \\
{\cal B} & \sim 1 \otimes_q \CB & \hookrightarrow   {\cal A} \otimes_q {\cal B}
\end{eqnarray}
we learn immediately that $q_A=p=q_B$. It remains now to
look for possible values of $q$ and $s$. The test of consistency is, 
as usually, the Leibniz rule.

We take the product $(a_1 \otimes_q b_1)(a_2 \otimes_q b_2)$ and
apply $d$ (defined as above) to it. This could be done
in two ways, first, we multiply the components and then apply $d$,
using the Leibniz rule for each of the algebra. The other way
is that we first apply $d$, using the above definition and then
multiply the result.

The table below shows the coefficients obtained in each way::

\medskip

\begin{tabular}{|r||c|c|c|c|}
\hline
 & $da_1 a_2 \otimes_q b_1 b_2$ & $a_1 da_2 \otimes_q b_1 b_2$ &
$a_1 a_2 \otimes_q db_1 b_2$ & $a_1 a_2 \otimes_q b_1 db_2$ \\
\hline \hline
I  & $q^{|a_2| |b_1|}$ & $q^{|a_2| |b_1|} p^{|a_1|} $ &
$q^{|a_2| |b_1|} s^{|a_1|+|a_2|} $ &
$q^{|a_2| |b_1|} s^{|a_1|+|a_2|} p^{|b_1|}$ \\ \hline
II & $q^{|a_2| |b_1|}$ & $ p^{|a_1|+|b_1|} q^{|b_1|(|a_1|+1)}$ &
$q^{(|b_1|+1)|a_2|} s^{|a_1|}$ &
$p^{|a_1|+|b_1|} s^{|a_2|} q^{|b_1| |a_2|}$ \\ \hline
\end{tabular}

\medskip

Of course, the consistency requires that the values in the first row must 
coincide with these in the second row. Let us compare: from the second 
column we get that $p = q^{-1}$, the third determines $s=q$ and from 
the last one we obtain $p=s$. However, this is possible only in the
case $p=q=s=-1$!. Therefore, the procedure of tensoring the
differential algebras is natural only in the standard $d^2=0$
situation.

\section{An example - but not a counterexample.}

Let us see a simple example illustrating the problem. We take
$N=3$, then we have two primitive roots of unity $j$ and 
$j^2=\bar{j}$. 

Let us take the simplest (but nontrivial) $d^3 \equiv 0$ $j$-differential
algebra on the real line. Our starting point is the standard definition of $d$
on the algebra of $C^\infty(\R)$. For $f \in C^\infty(\R)$ we define:
\begin{equation}
df \; = \; dx f',
\end{equation}
where $f'$ is the usual derivative and $dx$ is the one-form that
generates the bimodule of one-forms. Applying $d$ to both sides we get:
\begin{equation}
d^2f \; = \; d^2x f' + j dx dx f''.
\end{equation}

where we have assumed the $j$-Leibniz rule. The forms $d^2x$ and $dx dx$ 
are two independent generators of $\Omega^2_j(\R) $. Consider now the 
general $j$-Leibniz rule. First, for the product of two functions
$f$ and $g$ we have:
$$ d(fg) = dx (fg)' = dx (f'g + fg'),$$
however, on the other hand, from the $j$-Leibniz rule:
$$ d(fg) = (df) g + f (dg) = dx f' g + f dx g'.$$
This gives us the first rule for the left and right multiplication of
differential forms:
\begin{equation}
f dx \; = \; dx f, \label{bimo1}
\end{equation}
which is the same as for the standard calculus. However, if we apply
$d$ to the both sides of (\ref{bimo1}) we get:
\begin{equation}
 dx dx f' + f d^2x \; = \; d^2 x f + j dx dx f',
\end{equation}
where we have already used (\ref{bimo1}) when appropriate. Finally, we
obtain the second relation:
\begin{equation}
f d^2x \; = \; d^2x f + (j-1) dx dx f'. \label{bimo2}
\end{equation}

To obtain some further relations between the differential forms, we
apply $d$ for the third time, to both sides of (\ref{bimo2}), then
taking into account that $d^3 \equiv 0$ we see:
\begin{equation}
dx f' d^2x = j^2 d^2x dx f' + (j-1) d^2x dx f' +
j(j-1) dx d^2x + j^2(j-1) dx dx dx f''.
\end{equation}
Again, using both (\ref{bimo1}) and (\ref{bimo2}) we get:
\begin{equation}
-3 dx dx dx f'' -2 j^2 (dx d^2x -j d^2x dx) f' = 0,
\end{equation}
from which we immediately get:
\begin{eqnarray}
& dx dx dx  =  0, \\
& dx d^2x  =  j d^2x dx.
\end{eqnarray}

These are the only rules that are {\em necessary} for the construction
of a consistent $j$-differential algebra on $\R$. The obtained differential 
complex is infinite-dimensional, as there are no restrictions on products 
of $d^2 x$. However, without breaking the consistency, one may  impose more
rules, for instance, one may
postulate that $d^2 x d^2x = 0$ or $d^2x d^2x d^2x =0$. 

The one-form $dx$ is hermitian, $dx^* = dx$ whereas $d^2x$ 
is antihermitian: $(d^2x)^* = j^2 d^2x$.

Of course, the whole procedure could be carried out with $\bar{j}$
instead of $j$, giving us a different $d^3=0$ differential structure on 
$\R$. 

So far, we have been able to construct the $j$- (or $\bar{j}$-) differential
structures on the real line. The natural question is then, how can one
extend this construction to the plane or, more generally, $\R^n$. As
we have seen in the previous section, there is no canonical way to do it.
This does not mean, however, that such structure does not exist. 

Below, we shall provide an example of the $j$-differential structure on
the plane, which is built out of $j$ and $\bar{j}$ differential structures
on the line. This will not be, however, a counterexample to our {\em no-go}
lemma.

Consider the algebra of functions on the plane, and the differential
calculus defined similarly as in the one-dimensional case:
\begin{eqnarray}
df & = & dx f_x + dy f_y, \label{rule1} \\
d^2f & = & d^2x f_x + d^2y f_y + j (dx)^2 f_{xx} +
+j (dy)^2 f_{yy} + j (dx\, dy + dy\, dx) f_{xy}
\end{eqnarray}

We notice here that the $j$-differential algebras on the real line
are embedded in a natural way in $\Omega_j(\R^2)$ as constructed above.
However, we encounter problems when we come to the commutation relations 
between $x,y$ and their derivatives. From Eq.(\ref{rule1}) we have
\begin{eqnarray}
x\, dy = dy\, x, & \null &  y\, dx = dx\, y.
\end{eqnarray}

If we differentiate it (using the $j$-Leibniz rule) we obtain
first set of the restrictions:

\begin{eqnarray}
dx\, dy - j dy\, dx & = & d^2y\, x - x\, d^2y, \\
dy\, dx - j dx\, dy & = & d^2x\, y - y\, d^2x.
\end{eqnarray}

We immediately see that at least in one row the left-hand side
does not vanish and therefore either $x$ does not commute with
$d^2y$ or $y$ with $d^2x$. This proves that such structure does 
not come from the anyonic tensor product, even though it is 
well-defined and contains $j$-differential algebras on the real line.

\section{Conclusions}

What we have shown in this Letter has some profound consequences 
for the theory of $q$-differential algebras. First, we have no natural
procedure of tensoring such structures, while keeping all the
{\em natural} rules. Of course, one could attempt to relax one
or another of them, however, there is still ambiguity, which
should be kept and which could be allowed to change.

Having no canonical construction for the $q$-differential structures
on the products of spaces is a major drawback of the theory. 

Of course, let us notice that the problem arises only if we 
take into account the algebra structure, as for the 
simple $q$-generalizations of cochains and cohomology \cite{Kap}
the construction of the tensor product might be possible.

The world of $q$-deformations, and $q$-differential structures seems, at
first sight, similar to the well-known commutative and anticommutative
objects. What we see at least in the case of the $q$-differential 
algebras is that $q=-1$ case is very special.

{\bf Acknowledgements:} It is a pleasure to thank R.Kerner for 
discussions and M.Dubois-Violette for remarks on this topic.

\end{document}